\documentclass{article}
\usepackage[english]{babel}

\usepackage{amsmath,amssymb,setspace}
\usepackage{amsthm}
\usepackage[english]{babel}
\usepackage{amsfonts}
\usepackage{calligra}
\usepackage[T1]{fontenc}
\usepackage{xstring}

\usepackage[colorlinks,citecolor=blue,urlcolor=blue,bookmarks=false,hypertexnames=true]{hyperref}
\hypersetup{
    colorlinks=true,
    linkcolor=blue,
    filecolor=magenta,
    urlcolor=cyan,
}
\usepackage{graphics}
\usepackage[numbers,sort&compress]{natbib}
\usepackage{subcaption}
\usepackage{multicol}
\date{}
\usepackage{longtable}
\usepackage{float}
\usepackage[margin=0.7 in]{geometry}
\usepackage{graphicx}
\usepackage{mathtools}

\newtheorem{definition}{Definition}
\theoremstyle{plain}
\theoremstyle{definition}
\theoremstyle{remark}
\newtheorem{theorem}{Theorem}

\newtheorem{remark}{Remark}
\newtheorem{example}{Example}

\title{Fuzzification of  Fractal Calculus}
\author{Alireza Khalili Golmankhaneh \footnote{Corresponding author} \\
Department of Physics, Urmia Branch, Islamic Azad University, Urmia 63896,  Iran\\
	Kerri Welch\\
Faculty at California Institute of Integral Studies, San Francisco, CA, U.S.A.\\
Cristina Serpa\\
Instituto Superior de Engenharia de Lisboa (ISEL), \\Instituto Polit\'{e}cnico de Lisboa, Lisbon, Portugal\\
Centro de Matem\'{a}tica, Aplica\c{c}\~{o}es Fundamentais e Investiga\c{c}\~{a}o Operacional (CMAFcIO),\\ Faculdade de Ci\^{e}ncias, Universidade de Lisboa, Lisbon, Portugal\\
Palle E. T. J{\o}rgensen\\
Department of Mathematics, The University of Iowa, Iowa City, IA
52242-1419, U.S.A.
}

\date{\today}
\begin{document}

\maketitle

\begin{abstract}
In this manuscript, fractal and fuzzy calculus are summarized.  Fuzzy calculus in terms of fractal limit, continuity, its derivative, and integral are formulated. The fractal fuzzy calculus is a new framework that includes fractal fuzzy derivatives and fractal fuzzy integral. In this framework, fuzzy number-valued functions with fractal support are the solutions of fractal fuzzy differential equations. Different kinds of fractal fuzzy differential equations are given and solved.
\end{abstract}
\section{Introduction}
Fractal geometry mathematically describes complex shapes that are not described by Euclidean geometry \cite{b-1}. These shapes are found in nature such as clouds, mountains, lightning and etc. which are called fractals \cite{b-6}. The most important properties of fractals are self-similarity, and have non-integer dimensions. Fractals are non-differentiable in the sense of ordinary calculus since they have a rough structure rather smooth. Their fractal dimensions exceed their topological dimensions and appear similar at various scales \cite{falconer1999techniques,ma-12}. Fractals have different measures like Hausdorff's measure. In this context,  ordinary calculus which is based on length, area, and volume fails to define derivatives and integrals on them \cite{b-2}.\\
 Many researchers have tried to formulate analysis on  fractals in order to explain their physical properties \cite{samayoa2020fractal}, i.e. harmonic analysis \cite{ma-8,ma-13,freiberg2002harmonic}, measure theory \cite{ma-3}, fractional Brownian motion and probability-theoretical approaches \cite{ma-7,lee2022propagation}, fractional space \cite{stillinger1977axiomatic}, fractional calculus \cite{ma-6,trifce2020fractional}.\\
In seminal papers, ordinary calculus was adopted to define their  derivatives and integrals of  functions with fractal support, like Cantor sets and Koch curves \cite{parvate2009calculus,AD-2,parvate2011calculus}. This new framework which is a generalization of ordinary calculus is called fractal calculus or $F^{\alpha}$-calculus. Fractal calculus is  simple, constructive, and algorithmic and applied in physics \cite{Alireza-book}.
Fractal calculus  was developed in different branches such as stability of solutions of fractal differential equations, nonlocal reverse Minkowski's fractal integral inequalities, and properties of staircase function
\cite{khalili2021hyers,khalili2019fractalcat,shapovalov2021invariance,rahman2021nonlocal,cetinkaya2021general,wibowo2021relationship}.
 Fractal calculus was generalized on fractal cubes and tartan Cantor spaces and Laplace equations on fractal cubes  were solved  \cite{golmankhaneh2018fractalt,khalili2021laplace}. Fractal derivatives and integrals were worked out fractal interpolation functions and Weierstrass functions \cite{gowrisankar2021fractal}. Random variables, stochastic process and stable distributions on fractal were defined, and corresponding stochastic differential equations were solved
\cite{khalili2019random,golmankhaneh2020stochastic}. Fractal Laplace, Fourier and Sumudu transforms were defined in order to solve fractal differential equation and applied in electrical circuits, economy and in dynamics \cite{Rewid3,banchuin2022noise,golmankhaneh2016non,golmankhaneh2019sumudu,Fourier1,banchuin20224noise,khalili2021economic}.

 Fractal anomalous diffusion has been formulated as a diffusion process in fractal media and it has a power law relationship between the mean squared displacement and time \cite{Alireza-Fernandez-1,golmankhaneh2018sub,golmankhaneh2021equilibrium}.\\
 Fuzzy sets, numbers, fuzzy-valued functions, fuzzy derivatives, and integrals  were introduced and applied to model the processes with uncertainty in science, physical science, engineering, and social science   \cite{kloeden1994metric,allahviranloo2021fuzzy,anastassiou2010fuzzy,lakshmikantham2004theory}.
A linear second-order differential equation with constant  coefficients with boundary values expressed by fuzzy numbers have been solved \cite{gasilov2014solution}. The fuzzy optimal control problem has been considered to optimize the expected values of the appropriate objective fuzzy functions \cite{zarei2022suboptimal}. The differentiability of fuzzy number-valued functions based on the Hausdorff distance between fuzzy numbers has been suggested  \cite{khastan2022new}. First order linear fuzzy differential equations under differential inclusions and strongly generalized differentiability approaches have been studied \cite{khastan2020linear}. Linear fuzzy differential equations appying the concept of generalized differentiability and  conditions for the existence of solutions have been investigated \cite{khastan2016solutions}.

 A First-order fuzzy differential equation (FDE) with fuzzy initial value was solved \cite{sedaghatfar5method}.\\
In this paper, we introduce a new framework which is a generalization of fractal calculus to include  fuzzy-valued functions. \\
The plan of the paper is as follows:\\
In Section \ref{1g}, we summarize the fractal calculus and fuzzy calculus. Fractal fuzzy calculus is formulated and defined in Section \ref{2g}. In Section \ref{3g}, $\alpha$-order fractal fuzzy differential equations are suggested  and solved. Section \ref{4g} is devoted to conclusion.

\section{Preliminaries \label{1g}}
In this section we  summarize the fractal calculus on fractal curves \cite{parvate2009calculus,AD-2,parvate2011calculus,Alireza-book}.

\subsection{Fractal calculus on fractal curves}
\begin{definition}
For a fractal curve $F$ and a subdivision $P_{[a,b]}, [a,b]\in [a_{0},b_{0}] \in \mathbb{R} $, the mass function is defined by
\begin{equation}
\gamma^{\alpha}(F,a,b)=\lim_{\delta\rightarrow0} \inf_{|P|\leq \delta}\sum_{i=0}^{n-1}
\frac{|\mathbf{w}(t_{i+1})-\mathbf{w}(t_{i})|^{\alpha}}{\Gamma(\alpha+1)},
\end{equation}
where $|.|$ denotes the Euclidean norm on $\mathbb{R}^{n}$, $1\leq \alpha\leq n$, $P_{[a,b]}=\{a=t_{0},...,t_{n}=b\}$, and $|P|=\max_{0\leq i\leq n-1}(t_{i+1}-t_{i})$ for a subdivision $P$.
\end{definition}
\begin{definition}
The $\gamma$-dimension of $F$ is defined by
\begin{align}
  \dim_{\gamma}(F)&=\inf\{\alpha:\gamma^{\alpha}(F,a,b)=0\}\nonumber\\&=
\sup\{\alpha:\gamma^{\alpha}(F,a,b)=\infty\}
\end{align}
\end{definition}

\begin{definition}
The rise function of a fractal curve $F$ is defined by
\begin{equation}
  S_{F}^{\alpha}(u)=\left\{
                      \begin{array}{ll}
                        \gamma^{\alpha}(F,p_{0},u), & u\geq p_{0} ; \\
                        -\gamma^{\alpha}(F,u,p_{0}), & u<p_{0}.
                      \end{array}
                    \right.
\end{equation}
where $u\in [a_{0},b_{0}]$, and $S_{F}^{\alpha}(u)$  gives the mass of the fractal curve $F$ upto point $u$.
\end{definition}
\begin{definition}
Let be a function $f:F\rightarrow \mathbb{R}$. Then  $F$-limit of $f$ as $\theta'\rightarrow \theta$  through points of $F$ is $l$, if for given $\epsilon$ there exists $\delta>$ such that
\begin{equation}
  \theta'\in F~~and~~|\theta'-\theta|<\delta\Rightarrow |f(\theta')-l|<\epsilon
\end{equation}
or
\begin{equation}
  \underset{ \theta'\rightarrow \theta}{F_{-}lim}f(\theta')=l.
\end{equation}
\end{definition}
\begin{definition}
A function $f:F\rightarrow \mathbb{R}$  is said to be $F$-continuous at $\theta$ if
\begin{equation}
  \underset{ \theta'\rightarrow \theta}{F_{-}lim}f(\theta')=f(\theta).
\end{equation}
\end{definition}
\begin{definition}
 The fractal derivative $F^{\alpha}$-derivative is defined by
\begin{equation}
  D_{F}^{\alpha}f(\theta)=\underset{ \theta'\rightarrow \theta}{F_{-}lim}~
\frac{f(\theta')-f(\theta)}{J(\theta')-J(\theta)},
\end{equation}
where $F_{-}lim$ indicates the fractal limit (see in \cite{parvate2011calculus}), $\mathbf{w}(u)=\theta$ and $S_{F}^{\alpha}(u)=J(\theta)$.
\end{definition}
\begin{remark}
We note that the Euclidean distance from origin upto a point $\theta=\mathbf{w}(u)$ is given by $L(\theta)=L(\mathbf{w}(u))=|\mathbf{w}(u)|.$
\end{remark}
\begin{definition}
The fractal integral or $F^{\alpha}$-integral is defined by
\begin{align}
  \int_{C(a,b)}f(\theta)d_{F}^{\alpha}\theta&=\sup_{P[a,b]}\sum_{i=0}^{n-1}
\inf_{\theta\in C(t_{i},t_{i+1})}f(\theta)(J(\theta_{i+1})-J(\theta_{i}))\nonumber\\&=
\inf_{P[a,b]}\sum_{i=0}^{n-1}
\sup_{\theta\in C(t_{i},t_{i+1})}f(\theta)(J(\theta_{i+1})-J(\theta_{i})),
\end{align}
where $t_{i}=\mathbf{w}^{-1}(\theta_{i})$, and $C(a,b)$ is the section of the curve lying between points $\mathbf{w}(a)$ and $\mathbf{w}(b)$ on the fractal curve $F$ \cite{parvate2011calculus}.
\end{definition}

\subsection{Fuzzy calculus on real-line}
In this section, we review fuzzy calculus which will be used to fuzzification of the fractal calculus \cite{kloeden1994metric,allahviranloo2021fuzzy,anastassiou2010fuzzy,lakshmikantham2004theory}.\\
A generalized  Hukuhara difference for fuzzy sets and a new generalized differentiability concepts for fuzzy valued functions were given in  \cite{bede2013generalized,stefanini2009generalized}.

\begin{definition}
Let $X\neq\emptyset$. Then, a set $A\subset X$ is characterized by its membership function $u_{A}(x):X\rightarrow [0,1]$. Thus $u_{A}(x)$ is the degree of membership of element $x$ in the fuzzy set $A$ for each $x\in X$.
\end{definition}
\begin{definition}
Let $A$ be a fuzzy subset of a real number $u_{A}(x):\mathbb{R}\rightarrow [0,1]$. Then
$A$ is called fuzzy number if it satisfies the following axioms
\begin{enumerate}
  \item $A$ is  normal. It means that there exists $x_{0}$ in $ \mathbb{R}$, such that $u_{A}(x_{0})=1$.
  \item $A$ is convex, namely,
\begin{equation}
  u_{A}(tx+(1-t)y)\geq \min \{u_{A}(x),u_{A}(y)\},~\forall t\in [0,1],~ x,~ y\in \mathbb{R}.
\end{equation}
  \item $u_{A}(x)$ is upper semi continuous on $\mathbb{R}$, viz, for given $\epsilon>0$, there exists $\delta>0$ such that
\begin{equation}
  |x-x_{0}|<\delta \Rightarrow u_{A}(x)-u_{A}(x_{0})<\epsilon.
\end{equation}
\item The support of $u_{A}(x)$ is compact. e.g.
\begin{equation}
  supp(u_{A}(x))=cl_{\mathbb{R}}\{x\in \mathbb{R};u_{A}(x)> 0\}
\end{equation}
is compact.
\end{enumerate}
\end{definition}
\begin{definition}\label{oo}
A fuzzy number $A$ is determined by a pair of functions $A=(A^{-}(r),A^{+}(r))$, with $A^{-}(r),A^{+}(r):[0,1]\rightarrow \mathbb{R}$ that satisfies the following condition:
\begin{enumerate}
  \item $A^{-}(r)=A_{r}^{-}\in \mathbb{R}$ is a bounded, monotonic, increasing, left continuous function in $(0,1]$ and it is right-continuous at $0$.
  \item $A^{+}(r)=A_{r}^{+}\in \mathbb{R}$ is a bounded, monotonic, decreasing, left continuous function in $(0,1]$ and it is right-continuous at $0$.
  \item  For $r\in (0,1]$ we have $A^{-}(r)\leq A^{+}(r)$.
\end{enumerate}
The Definition \ref{oo} is called parametric form of fuzzy numbers.
\end{definition}

\begin{definition}
The $r$-cut of fuzzy number $A$ is defined and called level wise form by
\begin{equation}
  [A]_{r}=A_{r}=\{x\in \mathbb{R};u_{A}(x)\geq r\}
\end{equation}
where $A_{r}$ is closed interval $A_{r}=[A_{r}^{-},A_{r}^{+}]$ for any $r\in[0,1]$. We note that  $[A]_{0}=supp(u_{A}(x))$ and $F_{\mathbb{R}}$ denote space of fuzzy number.
\end{definition}

\begin{definition}\label{yyhbvgt7}
For every $A,B \in F_{\mathbb{R}}$ and $\lambda\in \mathbb{R},~r\in [0,1]$ the addition and scalar multiplication is defined by
\begin{equation}
  (A\oplus B)_{r}=A_{r}+B_{r},~~~(\lambda\odot A)_{r}=\lambda A_{r}.
\end{equation}
\end{definition}
\begin{definition}
The Hausdorff distance between two fuzzy numbers $A,~B$ using their $r$-cuts is defined by
\begin{equation}
 d_{H}(A,B)=\sup_{0\leq r\leq 1} \max\{|A_{r}^{-}-B_{r}^{-}|,|A_{r}^{+}-B_{r}^{+}|\}
\end{equation}
\end{definition}
\begin{remark}
  The set of fuzzy numbers $(F_{\mathbb{R}},d)$ with addition and scalar multiplication given in Definition \ref{yyhbvgt7}, is a complete metric space.
\end{remark}
\begin{definition}
Consider $A,~B\in F_{\mathbb{R}}$. The Hukuhara difference of $A,B$ is defined by
\begin{equation}
  C=A\ominus B,
\end{equation}
if $A=B\oplus C$.
\end{definition}
\begin{definition}
 Consider a fuzzy number valued function $f:\mathbb{R}\rightarrow F_{\mathbb{R}}$ and $x_{0}\in \mathbb{R}$, then $l$ is called limit of $f$ at point $x_{0}$ if for every given $\epsilon>0$, there exist $\delta>0$ such that \cite{anastassiou2010fuzzy,ghaffari2022generalized,allahviranloo2021fuzzy}
\begin{equation}
   0<|x-x_{0}|<\delta\Rightarrow d_{H}(f(x),l)<\epsilon,
\end{equation}
or,
\begin{equation}
  \lim_{x\rightarrow x_{0}} f(x)=l,
\end{equation}
if it exists, where $d_{H}$ is the Hausdorff distance.
\end{definition}

\begin{definition}
The fuzzy function $f$ is called fuzzy continuous if \cite{anastassiou2010fuzzy,ghaffari2022generalized}
\begin{equation}
  \lim_{x\rightarrow x_{0}} f(x)=f(x_{0}).
\end{equation}
\end{definition}
\begin{definition}\label{km8558uyh}
A fuzzy number valued function $f:\mathbb{R}\rightarrow F_{\mathbb{R}}$ is called Hukuhara differentiable if there exist $f'(x)\in F_{\mathbb{R}}$ such that
\begin{itemize}
  \item Case 1.($I$-differentiable)
 \begin{equation}
  f'(x)=\lim_{y\rightarrow x} \frac{f(y)\ominus f(x)}{y-x},~~~y>x
\end{equation}
  \item Case 2.($II$-differentiable)
\begin{equation}
 f'(x)=\lim_{y\rightarrow x} \frac{f(x)\ominus f(y)}{y-x},~~~y>x
\end{equation}
where $f'(x)$ is called the fuzzy derivative of $f$ at $x$.
\end{itemize}

\end{definition}

\begin{theorem}
 Let a fuzzy number valued function $f:\mathbb{R}\rightarrow F_{\mathbb{R}}$ be denoted by $f(x)=(\underline{f}(x,r),\overline{f}(x,r))$ for each $r\in [0,1]$ \cite{chalco2008new}. Then \cite{khastan2011variation,khastan2016solutions}
\begin{enumerate}
  \item If $f$ is $I$-differentiable, then we have
\begin{equation}
  f'(x)=(\underline{f}'(x,r),\overline{f}'(x,r)).
\end{equation}
  \item If $f$ is $II$-differentiable, then we have
\begin{equation}
  f'(x)=(\overline{f}'(x,r),\underline{f}'(x,r)).
\end{equation}
\end{enumerate}
\end{theorem}

\begin{definition}
Let $f(x)$ be a fuzzy number-valued function. Then the fuzzy Riemann integral is defined as \cite{allahviranloo2021fuzzy}
\begin{equation}
  J=FR\int_{a}^{b}f(x)dx=\oplus\sum_{i=0}^{n}\Delta x_{i}\odot f(x_{i}),
\end{equation}
 where $\Delta x_{i}=x_{i+1}-x_{i}$ and $\{a=x_{0}<x_{1}<...< x_{n}=b\}$ is a partition of $I=[a,b]$. \\
The fuzzy Riemann integral of $f(x)$ is $J$ if for every given $\epsilon>0$, there exist $\delta>0$ such as
\begin{equation}
  d_{H}\bigg(\oplus\sum_{i=0}^{n}\Delta x_{i}\odot f(x_{i}),J\bigg)<\epsilon.
\end{equation}
where $J$ is a fuzzy number.
\end{definition}

\begin{definition}
Let $f:I\rightarrow F_{\mathbb{R}}$ be a triangular number-valued function and $f(x)=(f_{1}(x),f_{2}(x),f_{3}(x))$ and $x_{0}\in I$. Then the fuzzy integral is defined by \cite{ghaffari2022generalized}
\begin{equation}
  \int_{a}^{b}f(x)dx=\bigg(\int_{a}^{b}f_{1}(x)dx,\int_{a}^{b}f_{2}(x)dx,\int_{a}^{b}f_{3}(x)dx\bigg)
\end{equation}
\end{definition}

\section{Fuzzy fractal calculus on fractal curves \label{2g}}
 In this section, we introduce fractal fuzzy calculus.
\begin{definition}
Let $f(\theta):F\rightarrow F_{\mathbb{R}}$ be a number-valued function on a fractal curve $F$. Then the fuzzy $F$-limit of $f$ at $\theta_{0}$ through $F$ is $l$, if for a given $\epsilon>0$, there exist $\delta>0$, such that
\begin{equation}
  \theta\in F,~~~and~~~|\theta-\theta_{0}|<\delta\Rightarrow d_{H}(f(\theta),l)<\epsilon,
\end{equation}
or
\begin{equation}
   \underset{ \theta\rightarrow \theta_{0}}{FF_{-}lim}f(\theta')=l
\end{equation}
where $d_{H}$ is the Hausdorff distance.
\end{definition}

\begin{definition}
Let  $f(\theta):F\rightarrow F_{\mathbb{R}}$ be a number-valued function on $F$. Then, $f$ is called fuzzy $F$-continuous if
\begin{equation}
  \underset{ \theta\rightarrow \theta_{0}}{FF_{-}lim}f(\theta)=f(\theta_{0}).
\end{equation}
\end{definition}
\begin{definition}\label{Fuzzyfractalderivative}
Let $f(\theta):F\rightarrow F_{\mathbb{R}}$ be a number-valued function. Then fractal Hukuhara differentive at $\theta_{0}\in F$ is defined by
\begin{itemize}
  \item Case 1. ($I$-$F^{\alpha}$-differentiable)
\begin{equation}
  D_{F,H}^{\alpha}f(\theta_{0})=\underset{ \theta\rightarrow \theta_{0}}{FF_{-}lim}~ \frac{f(\theta)\ominus f(\theta_{0})}{J(\theta)-J(\theta_{0})},~~~\theta>\theta_{0}.
\end{equation}
  \item Case 2. ($II$-$F^{\alpha}$-differentiable)
\begin{equation}
  D_{F,H}^{\alpha}f(\theta_{0})=\underset{ \theta\rightarrow \theta_{0}}{FF_{-}lim}~ \frac{f(\theta_{0})\ominus f(\theta)}{J(\theta)-J(\theta_{0})},~~~~\theta>\theta_{0}.
\end{equation}
where $D_{F,H}^{\alpha}f(\theta_{0})$ is a fuzzy number.
\end{itemize}
\end{definition}

\begin{definition}
Let $f(x)$ be a fractal fuzzy number valued function. Then the fractal fuzzy Riemann integral is defined as \cite{allahviranloo2021fuzzy}
\begin{equation}
  J=FFR\int_{C(a,b)}f(\theta)d_{F}^{\alpha}\theta=\oplus\sum_{i=0}^{n}\Delta J_{i}\odot f(\theta_{i}),
\end{equation}
  The fractal fuzzy Riemann integral of $f(\theta)$ is $J$ if for a given $\epsilon>0$, there exist $\delta>0$ such as
\begin{equation}
  d_{H}\bigg(\oplus\sum_{i=0}^{n}\Delta J_{i}\odot f(\theta_{i}),J\bigg)<\epsilon.
\end{equation}
\end{definition}
\begin{definition}
Let $f:F\rightarrow F_{\mathbb{R}}$ be a fractal triangular number-valued function, $f(\theta)=(f_{1}(\theta),f_{2}(\theta),f_{3}(\theta))$, and $\theta_{0}\in F$, then
\begin{equation}
  \int_{C(a,b)}f(\theta)d_{F}^{\alpha}\theta=
\bigg(\int_{C(a,b)}f_{1}(\theta)d_{F}^{\alpha}\theta,
\int_{C(a,b)}f_{2}(\theta)d_{F}^{\alpha}\theta,
\int_{C(a,b)}f_{3}(\theta)d_{F}^{\alpha}\theta\bigg).
\end{equation}
\end{definition}

\section{Fractal fuzzy differential equations \label{3g}}
First order linear fuzzy differential equations by using the generalized
differentiability concept were solved  \cite{khastan2011variation,khastan2016solutions}.
In this section, a $\alpha$-order fuzzy differential equation (F.D.E) is given. Then it is changed by its equivalent parametric form, and a new system, which contains two fractal differential equations, is solved.\\
Consider the following fractal fuzzy differential equation with  initial condition:
\begin{equation}\label{yyh85}
  D_{F,H}^{\alpha}x(\theta)=f(J(\theta),x(\theta)),~~~\tilde{x}(\theta_{0})=
\tilde{x}_{0},~~~\theta\in F,
\end{equation}
where
$f:F\times F_{\mathbb{R}}\rightarrow F_{\mathbb{R}}$ is a fuzzy-valued function and $\tilde{x}_{0}\in F_{\mathbb{R}}$. To solve Eq.\eqref{yyh85},  first we solve $1$-cut and $0$-cut of Eq.\eqref{yyh85} as the following form
\begin{equation}\label{gggsaws1}
\left\{
  \begin{array}{ll}
    (D_{F,H}^{\alpha}x)^{[1]}(\theta)=f^{[1]}(J(\theta),x(\theta)), &  \\
&\\
    x^{[1]}(\theta_{0})=\tilde{x}_{0}^{[1]}~~\theta_{0}\in[0,\Theta]. &
  \end{array}
\right.
\end{equation}
and
\begin{equation}\label{gggsaws2}
\left\{
  \begin{array}{ll}
    (D_{F,H}^{\alpha}x)^{[0]}(\theta)=f^{[0]}(J(\theta),x(\theta)), &  \\
&\\
    x^{[0]}(\theta_{0})=\tilde{x}_{0}^{[0]}~~\theta_{0}\in[0,\Theta]. &
  \end{array}
\right.
\end{equation}
Then by solving Eqs.\eqref{gggsaws1} and \eqref{gggsaws2}, we can find $\tilde{x}(\theta)$ which is the solution of the fractal fuzzy differential equation Eq.\eqref{yyh85}.\\
Here we consider two cases:\\
Case (I): Suppose that $\tilde{x}(\theta)$ is $I$-$F^{\alpha}$-differentiable. Then, we can write
\begin{equation}\label{dds521}
  D_{F,H}^{\alpha}x(\theta)=[D_{F,H}^{\alpha}
\underline{x}(\theta,r),D_{F,H}^{\alpha}\overline{x}(\theta,r)].
\end{equation}
 In view of Eqs.\eqref{yyh85} and \eqref{dds521} for $r\in [0,1]$, we have
\begin{equation}
  \left\{
    \begin{array}{ll}
      D_{F,H}^{\alpha}
\underline{x}(\theta,r)=\underline{f}(\theta,r) & \theta_{0}\leq\theta \leq \Theta \\
&\\
      D_{F,H}^{\alpha}\overline{x}(\theta,r)=\overline{f}(\theta,r) & \theta_{0}\leq\theta \leq \Theta.
    \end{array}
  \right.
\end{equation}
Hence
\begin{equation}
\left\{
  \begin{array}{ll}
    [(1-r)
D_{F,H}^{\alpha}\underline{x^{[0]}}(\theta)+r D_{F,H}^{\alpha}\underline{x^{[1]}}(\theta)= (1-r)\underline{f^{[0]}}(\theta)+r(\underline{f^{[1]}})(\theta), & \\
& \\
   (1-r)
D_{F,H}^{\alpha}\overline{x^{[0]}}(\theta)+r D_{F,H}^{\alpha}\overline{x^{[1]}}(\theta)= (1-r)\overline{f^{[0]}}(\theta)+r(\overline{f^{[1]}})(\theta),
& \\
& \\
\underline{x}(\theta_{0},r)=(1-r)\underline{x^{[0]}}(\theta_{0})+r \underline{x^{[1]}}(\theta_{0}),& \\
& \\
\overline{x}(\theta_{0},r)=(1-r)\overline{x^{[0]}}(\theta_{0})+r \overline{x^{[1]}}(\theta_{0}).
  \end{array}
\right.
\end{equation}
It follows that
\begin{equation}\label{AQWZ1}
  \left\{
    \begin{array}{ll}
      D_{F,H}^{\alpha}\underline{x^{[0]}}(\theta)=\underline{f^{[0]}}(\theta), &  \\
&  \\
      D_{F,H}^{\alpha}\overline{x^{[0]}}(\theta)=\overline{f^{[0]}}(\theta) &\\
&\\
\underline{x^{[0]}}(\theta_{0})=\underline{x_{0}}^{[0]}&\\
&\\
\overline{x^{[0]}}(\theta_{0})=\overline{x_{0}}^{[0]},
    \end{array}
  \right.
\end{equation}
and

\begin{equation}\label{AQWZ2}
  \left\{
    \begin{array}{ll}
      D_{F,H}^{\alpha}\underline{x^{[1]}}(\theta)=\underline{f^{[1]}}(\theta), &  \\
&  \\
      D_{F,H}^{\alpha}\overline{x^{[1]}}(\theta)=\overline{f^{[1]}}(\theta) &\\
&\\
\underline{x^{[1]}}(\theta_{0})=\underline{x_{0}}^{[1]}&\\
&\\
\overline{x^{[1]}}(\theta_{0})=\overline{x_{0}}^{[1]},
    \end{array}
  \right.
\end{equation}
One can find $\underline{x^{[0]}}(\theta),\overline{x^{[0]}}(\theta),
\underline{x^{[1]}}(\theta),\overline{x^{[1]}}(\theta)$ by solving Eqs.\eqref{AQWZ1} and \eqref{AQWZ2}. Therefore we obtain the solution of Eq.\eqref{yyh85} using $0$-cut and $1$-cut solutions as follows:

\begin{equation}\label{Fistr}
  \tilde{x}(\theta)=[\underline{x}(\theta,r),\overline{x}(\theta,r)]=[(1-r)
\underline{x^{[0]}}(\theta)+r \underline{x^{[1]}}(\theta),(1-r)
\overline{x^{[0]}}(\theta)+r\overline{x^{[1]}}(\theta)].
\end{equation}
Case (II). Let $\tilde{x}(\theta)$ be $II$-$F^{\alpha}$-differentiable. Then, we can write
\begin{equation}\label{rdfrw}
  D_{F,H}^{\alpha}x(\theta)=[D_{F,H}^{\alpha}
\overline{x}(\theta,r),D_{F,H}^{\alpha}\underline{x}(\theta,r)].
\end{equation}
Likewise, Case (I), we have
\begin{equation}\label{IIAQWZ1}
  \left\{
    \begin{array}{ll}
      D_{F,H}^{\alpha}\underline{x^{[0]}}(\theta)=\overline{f^{[0]}}(\theta), &  \\
&  \\
      D_{F,H}^{\alpha}\overline{x^{[0]}}(\theta)=\underline{f^{[0]}}(\theta) &\\
&\\
\underline{x^{[0]}}(\theta_{0})=\overline{x_{0}}^{[0]}&\\
&\\
\overline{x^{[0]}}(\theta_{0})=\underline{x_{0}}^{[0]},
    \end{array}
  \right.
\end{equation}
and

\begin{equation}\label{IIAQWZ2}
  \left\{
    \begin{array}{ll}
      D_{F,H}^{\alpha}\underline{x^{[1]}}(\theta)=\overline{f^{[1]}}(\theta), &  \\
&  \\
      D_{F,H}^{\alpha}\overline{x^{[1]}}(\theta)=\underline{f^{[1]}}(\theta) &\\
&\\
\underline{x^{[1]}}(\theta_{0})=\overline{x_{0}}^{[1]}&\\
&\\
\overline{x^{[1]}}(\theta_{0})=\underline{x_{0}}^{[1]},
    \end{array}
  \right.
\end{equation}
 By solving the ordinary fractal differential equations \eqref{IIAQWZ1}) and \eqref{IIAQWZ2}, one may obtain the solution of FDE \eqref{yyh85} which is $II$-$F^{\alpha}$-differentiable as
\begin{equation}\label{Mohham}
  \tilde{x}(\theta)=[\overline{x}(\theta,r),\underline{x}(\theta,r)].
\end{equation}

\begin{example}
  Consider the fractal fuzzy differential equation as
\begin{equation}\label{UUUU}
   D_{F,H}^{\alpha}x(\theta)=x(\theta)+\tilde{c},
\end{equation}
with the conditions
\begin{equation}
 x(0,r)=[r,2-r],~~\tilde{c}=[r-1,1-r], r\in[0,1].
\end{equation}
Here we consider two cases.\\
Case I. Let $  \tilde{x}(\theta)$ be $I$-$F^{\alpha}$-differentiable. Then by using Eqs.\eqref{AQWZ1} and \eqref{AQWZ2} we arrive at
\begin{equation}\label{IIAQWZ1iu}
  \left\{
    \begin{array}{ll}
      D_{F,H}^{\alpha}\underline{x^{[0]}}(\theta)=\underline{x^{[0]}}(\theta)-1 &  \\
&  \\
      D_{F,H}^{\alpha}\overline{x^{[0]}}(\theta)=\overline{x^{[0]}}(\theta)+1 &\\
&\\
\underline{x^{[0]}}(\theta_{0})=0&\\
&\\
\overline{x^{[0]}}(\theta_{0})=2
    \end{array}
  \right.
\end{equation}
and

\begin{equation}\label{IIAQWZ2pp}
  \left\{
    \begin{array}{ll}
      D_{F,H}^{\alpha}\underline{x^{[1]}}(\theta)=\underline{x^{[1]}}(\theta) &  \\
&  \\
      D_{F,H}^{\alpha}\overline{x^{[1]}}(\theta)=\overline{x^{[1]}}(\theta) &\\
&\\
\underline{x^{[1]}}(\theta_{0})=1&\\
&\\
\overline{x^{[1]}}(\theta_{0})=1
    \end{array}
  \right.
\end{equation}
 By solving Eqs.\eqref{IIAQWZ1iu} and \eqref{IIAQWZ2pp}, we obtain
\begin{align}\label{Man2}
 \underline{x^{[0]}}(\theta)&=-\exp(J(\theta))+1,~~~
\overline{x^{[0]}}(\theta)=3\exp(J(\theta))-1\nonumber\\
\underline{x^{[1]}}(\theta)&=\exp(J(\theta))+1,~~~~~~
\overline{x^{[1]}}(\theta)=\exp(J(\theta))
\end{align}
 By substituting Eq.\eqref{Man2}  into Eq.\eqref{Fistr}, we get
\begin{equation}
  \tilde{x}(\theta)=[\exp(J(\theta))(2r-1)-r+1,r-\exp(J(\theta))(2r-3)-1].
\end{equation}
Case II. Let $\tilde{x}(\theta)$ be $II$-$F^{\alpha}$-differentiable. Then, by utilizing Eqs.\eqref{IIAQWZ1} and \eqref{IIAQWZ2} we get
\begin{equation}\label{kooo1}
  \left\{
    \begin{array}{ll}
      D_{F,H}^{\alpha}\underline{x^{[0]}}(\theta)=\overline{x^{[0]}}(\theta)+1, &  \\
&  \\
      D_{F,H}^{\alpha}\overline{x^{[0]}}(\theta)=\underline{x^{[0]}}(\theta)-1 &\\
&\\
\underline{x^{[0]}}(\theta_{0})=0&\\
&\\
\overline{x^{[0]}}(\theta_{0})=2,
    \end{array}
  \right.
\end{equation}
and

\begin{equation}\label{kooo2}
  \left\{
    \begin{array}{ll}
      D_{F,H}^{\alpha}\underline{x^{[1]}}(\theta)=\overline{x^{[1]}}(\theta), &  \\
&  \\
      D_{F,H}^{\alpha}\overline{x^{[1]}}(\theta)=\underline{x^{[1]}}(\theta) &\\
&\\
\underline{x^{[1]}}(\theta_{0})=1&\\
&\\
\overline{x^{[1]}}(\theta_{0})=1,
    \end{array}
  \right.
\end{equation}
Solving Eqs.\eqref{kooo1} and \eqref{kooo2} gives
\begin{align}\label{Jann}
  \underline{x^{[0]}}(\theta)&=
\exp(J(\theta))-r+\frac{(2r-2)}{\exp(J(\theta))}+1,~~~\overline{x^{[0]}}(\theta)=r
+\exp(J(\theta))-\frac{(2r-2)}{\exp(J(\theta))}-1\nonumber\\
\underline{x^{[1]}}(\theta)&=\exp(J(\theta)),~~~\hspace{3.5cm} \overline{x^{[1]}}(\theta)=\exp(J(\theta)).
\end{align}
By substituting Eq.\eqref{Jann} into Eq.\eqref{Mohham}, we find  the solution of Eq.\eqref{UUUU} as follows
\begin{equation}\label{wwwq}
 \tilde{x}(\theta)=\bigg[\exp(J(\theta))-r+\frac{2r-2}{\exp(J(\theta))}+1,
r+\exp(J(\theta))-\frac{2r-2}{\exp(J(\theta))}-1\bigg].
\end{equation}
In Figure \ref{dd}, we have plotted Eq.\eqref{wwwq} for the case of $r=0.3$.
\begin{figure}
  \centering
  \includegraphics[scale=0.8]{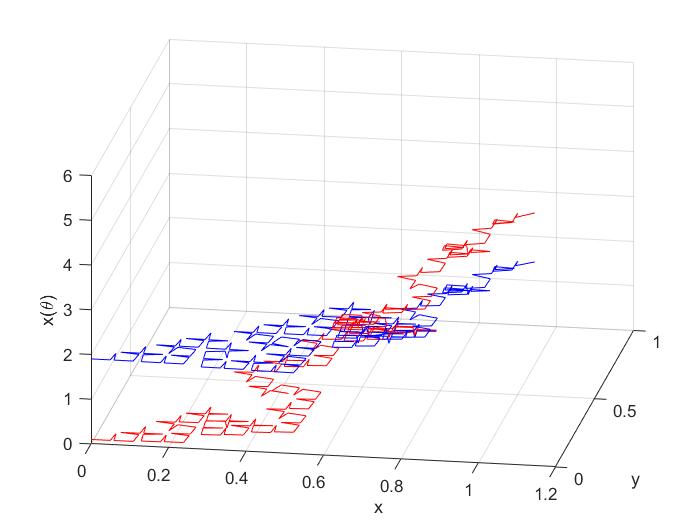}
  \caption{Graph of Eq.\eqref{wwwq} for $r=0.3$}\label{dd}
\end{figure}

\end{example}
\begin{example}
 Consider the fractal differential equation with fuzzy boundary condition:
\begin{equation}\label{Redree42587}
  \left\{
    \begin{array}{ll}
      (D_{F}^{\alpha})^{2}x(\theta)-4D_{F}^{\alpha}x(\theta)+4 x(\theta)=1-2 J^{2}(\theta), &  \\
   x(0)=(2,3,4), & \\
x(1)=(1,2,2.5),
    \end{array}
  \right.
\end{equation}
Let the solution be as $x(\theta)=x_{cr}(\theta)+\tilde{x}_{un}(\theta)$, where $x_{cr}$ is the crisp part of the solution and $\tilde{x}_{un}(\theta)$ is the uncertainty part \cite{gasilov2014solution}. First, we solve the following equation
\begin{equation}\label{EDxs}
  \left\{
    \begin{array}{ll}
      (D_{F}^{\alpha})^{2}x(\theta)-4D_{F}^{\alpha}x(\theta)+4 x(\theta)=1-2 J^{2}(\theta), &  \\
   x(0)=3, & \\
x(1)=2,
    \end{array}
  \right.
\end{equation}
The crisp  solution of Eq.\eqref{EDxs} is
\begin{equation}
  x_{cr}(\theta)=-\frac{1}{2}(J(\theta)+1)^{2}+3.5(1-J(\theta))\exp(2J(\theta))+
4J(\theta)\exp(2(J(\theta)-1))
\end{equation}
To determine the $ \tilde{x}_{un}(\theta)$, let us consider the fractal fuzzy homogeneous equations:
 \begin{equation}\label{Wsaq8}
  \left\{
    \begin{array}{ll}
      (D_{F}^{\alpha})^{2}x(\theta)-4D_{F}^{\alpha}x(\theta)+4 x(\theta)= 0, &  \\
   x(0)=(-1,0,1), & \\
x(1)=(-1,0,0.5)
    \end{array}
  \right.
\end{equation}
The linear independent solutions of Eq.\eqref{Wsaq8} are $x_{1}(\theta)=\exp(2J(\theta))$, and  $x_{2}(\theta)=J(\theta)\exp(2J(\theta))$
\end{example}
Then we have
\begin{align}
  &\mathbf{p}=(\exp(2J(\theta)),J(\theta)\exp(2J(\theta))),~~
\mathcal{M}=\left(
\begin{array}{cc}
1 & 0 \\
 e^{2} & e^{2} \\
\end{array}
\right),\nonumber\\&
~~~\mathbf{q}=\mathbf{p}\mathcal{M}^{-1}=((1-J(\theta))\exp(2J(\theta)),
J(\theta)\exp(2(J(\theta)-1)).
\end{align}
Thus we obtain
\begin{equation}
  \tilde{x}_{un}(\theta)=(1-J(\theta))\exp(2J(\theta))(-1,0,1)+
J(\theta)\exp(2(J(\theta)-1))(-1,0,0.5).
\end{equation}
To represent $ \tilde{x}_{un}(\theta)$ by $\kappa$-cuts, we can write
\begin{equation}\label{ttwqsrffd}
  \tilde{x}_{un,\kappa}(\theta)=(1-\kappa)[\underline{x_{un,0}}(\theta),
\overline{x_{un,0}}(\theta)].
\end{equation}
The $\kappa$-cuts representation solution of Eq.\eqref{Redree42587} is
\begin{equation}
 x_{\kappa}(\theta)=x_{cr}(\theta)+(1-\kappa)[\underline{x_{un,0}}(\theta),
\overline{x_{un,0}}(\theta)].
\end{equation}

\section{Conclusion \label{4g}}
In this paper, we have formulated the fractal fuzzy calculus which is a generalization of fractal calculus on fuzzy number-valued functions. Fractal calculus is a generalization of ordinary calculus which involves functions with a fractal domain such as Cantor set and Koch curves. Fractal fuzzy differential equations can be used to model uncertainty in the initial condition or dynamic of media with a fractal structure. The research in this direction is in progress.
\\

\textbf{Acknowledgements:} Cristina Serpa acknowledges partial funding by national funds through FCT-Foundation for Science and Technology, project reference: UIDB/04561/2020.
\bibliographystyle{elsarticle-num}
\bibliography{Fuzzy2}

\end{document}